\documentclass[12pt,a4paper]{article}
\usepackage{amsmath}

\usepackage{amscd,amsmath,amssymb,amsfonts}

\newtheorem{theorem}[subsection]{Theorem}
\newtheorem{definition}[subsection]{Definition}

\newtheorem{lemma}[subsection]{Lemma}
\newtheorem{corollary}[subsection]{Corollary}
\newtheorem{conjecture}[subsection]{Conjecture}
\newtheorem{remark}[subsection]{Remark}

\newenvironment{proof}{\bigskip\noindent\emph{Proof.}}{\bigbreak}

\title{Estimate of dimension of
Noether-Lefschetz locus for Beilinson-Hodge cycles
on open complete intersections}
\author{M. Asakura and S. Saito}
\date{}

\begin{document}
\maketitle
\tableofcontents
\medbreak
\begin{description}
\item[\S0]
Introduction
\item[\S1]
Jacobian rings for open complete intersections
\item[\S2]
Trivial part of cohomology of complements of hypersurfaces
\item[\S3]
Hodge theoretic implication of generalized Jacobian rings
\item[\S4]
Case of plane curves
\item[\S5]
Beilinson's Tate conjecture
\item[\S6]
Implication on injectivity of Chern class map for $K_1$ of surfaces
\item[]
References
\end{description}

\vskip 20pt

\renewcommand{\labelenumi}{(\theenumi)}


\def\indlim #1{\lim_{\longrightarrow \atop #1}}
\def\projlim #1{\lim_{\longleftarrow \atop #1}}
\def\qwith{\quad\hbox{with }}
\def\rmapo#1{\overset{#1}{\longrightarrow}}
\def\isom{\rmapo{\cong}}
\def\Spec{{\operatorname{Spec}}}
\def\Pic{{\mathrm{Pic}}}
\def\Ext{{\text{\mathrm{Ext}}}}
\def\NS{{\mathrm{NS}}}
\def\Picv{{\mathrm{Pic^0}}}
\def\Div{{\mathrm{Div}}}
\def\CH{{\mathrm{CH}}}
\def\deg{{\mathrm{deg}}}
\def\dim{{\operatorname{dim}}}
\def\codim{{\operatorname{codim}}}
\def\Coker{\mathrm{Coker}}
\def\dim{\mathrm{dim}}
\def\deg{\mathrm{deg}}
\def\Image{\mathrm{Image}}
\def\Ker{\mathrm{Ker}}
\def\Coker{\mathrm{Coker}}
\def\min{\mathrm{min}}
\def\Aut{\mathrm{Aut}}
\def\Hom{\mathrm{Hom}}
\def\Proj{\mathrm{Proj}}
\def\PGL{\mathrm{PGL}}
\def\Sym{\mathrm{Sym}}
\def\Gal{\mathrm{Gal}}
\def\GL{\mathrm{GL}}
\def\End{\mathrm{End}}
\def\Res#1{\mathrm{Res}_{#1}}
\def\ch{\mathrm{ch}}
\def\P{{\mathbb{P}}}
\def\bP{{\mathbb{P}}}
\def\C{{\mathbb{C}}}
\def\Q{{\mathbb{Q}}}
\def\bQ{{\mathbb{Q}}}
\def\Z{{\mathbb{Z}}}
\def\F{{\mathbb{F}}}
\def\bC{\mathbb C}
\def\bZ{\mathbb Z}
\def\cHom{{\mathcal{H}}om}
\def\cL{{\mathcal{L}}}
\def\L{{\mathcal{L}}}
\def\cE{{\mathcal{E}}}
\def\E{{\mathcal{E}}}
\def\cF{{\mathcal{F}}}
\def\cO{{\mathcal{O}}}
\def\cX{{\mathcal{X}}}
\def\cXB{{\mathcal{X}_B}}
\def\cC{{\mathcal{C}}}
\def\cZ{{\mathcal{Z}}}
\def\cZB{{\mathcal{Z}_B}}
\def\cZBj{{\mathcal{Z}_{B,j}}}
\def\cM{{\mathcal{M}}}
\def\cW{{\mathcal{W}}}
\def\cY{{\mathcal{Y}}}
\def\cV{{\mathcal{V}}}
\def\wtd#1{\widetilde{#1}}
\def\sE#1#2#3{E_{#1}^{#2,#3}}
\def\l{\ell}
\def\d{{\bold{d}}}
\def\e{{\bold{e}}}
\def\tH{\widetilde{H}}
\def\tR{\widetilde{R}}
\def\lra{\longrightarrow}
\def\ra{\rightarrow}
\def\la{\leftarrow}
\def\lla{\longleftarrow}
\def\Lra{\Longrightarrow}
\def\Lla{\Longleftarrow}
\def\da{\downarrow}
\def\hra{\hookrightarrow}
\def\lmt{\longmapsto}
\def\ot{\otimes}
\def\op{\oplus}
\def\wt#1{\widetilde{#1}}
\def\us#1#2{\underset{#1}{#2}}
\def\os#1#2{\overset{#1}{#2}}
\def\PE{\P(\cE)}
\def\qaq{\quad\hbox{ and }\quad}
\def\qfor{\quad\hbox{ for }}
\def\qif{\quad\hbox{ if }}
\def\ld{\lambda}
\def\scs{\; : \;}
\def\SS{\Sigma^*}
\def\ab{\underline{\alpha}}
\def\tB{\tilde{B}}
\def\tR{\tilde{R}}
\def\kX#1{\kappa_{\cX}^{#1}}
\def\tkX#1{\tilde{\kappa}_{\cX}^{#1}}
\def\tkY#1{\tilde{\kappa}_{\cY}^{#1}}
\def\kZ#1{\kappa_{\cZ^{(1)}}^{#1}}
\def\tkZ#1{\tilde{\kappa}_{\cZ^{(1)}}^{#1}}
\def\ddlog#1{\frac{d #1}{#1}}
\def\Sp{\Sigma'}
\def\Sb{\overline{\Sigma}}
\def\SSp{\Sigma^{'*}}
\def\SSb{\Sb^*}
\def\taub{\overline{\tau}}
\def\Pb{\overline{\mathbb P}}
\def\Rp#1#2{B'_{#1}(#2)}
\def\Rb#1#2{\overline{B}_{#1}(#2)}
\def\eb{\overline{\bold e}}
\def\pg#1{{\frak S}_{#1}}
\def\ul#1{\underline{#1}}
\def\ba{\bold{a}}
\def\bb{\bold{b}}
\def\ua{\underline{a}}
\def\ub{\underline{b}}
\def\be{\bold{e}}
\def\bd{\bold{d}}
\def\ue{\underline{e}}
\def\ud{\underline{d}}
\def\md{\delta_{min}}
\def\Om{\Omega}
\def\ld{\lambda}
\def\Wb{\overline{W}}
\def\XS{\cX_S}
\def\CS{\cC_S}
\def\ZS{\cZ_S}
\def\FS#1{F_S^{#1}}
\def\WS#1{\Omega_S^{#1}}
\def\GrFS#1{Gr_{F_S}^{#1}}
\def\otO{\otimes_{\cO}}
\def\otQ{\otimes{\mathbb Q}}
\def\cU{{\mathcal U}}
\def\onab{\overline{\nabla}}
\def\HO{H_{\cO}}
\def\HQ{H_{\mathbb Q}}
\def\HC{H_{\mathbb C}}
\def\HQcU#1#2{\HQ^{#1}(\cU/S)(#2)}
\def\HCcU#1{\HC^{#1}(\cU/S)}
\def\HOcU#1{\HO^{#1}(\cU/S)}
\def\HcU#1#2{H^{#1,#2}(\cU/S)}
\def\HU#1#2{H^{#1,#2}(U)}
\def\HccU#1#2{H_{\cO,c}^{#1,#2}(\cU/S)}
\def\HcZ#1#2{H^{#1,#2}(\cZ/S)}
\def\HcX#1#2{H^{#1,#2}(\cX/S)}
\def\HcX#1#2{H^{#1,#2}(\cX/S)}
\def\HcU#1#2{H^{#1,#2}(\cU/S)}
\def\HcXp#1#2{H^{#1,#2}(\cX/S)_{prim}}
\def\HXp#1#2{H^{#1,#2}(X)_{prim}}
\def\HcUan#1#2{H^{#1,#2}(\cU/S)^{an}}
\def\HQcU#1#2{\HQ^{#1}(\cU/S)(#2)}
\def\HQcX#1#2{\HQ^{#1}(\cX/S)(#2)}
\def\HQcU#1#2{\HQ^{#1}(\cU/S)(#2)}
\def\HQcXp#1#2{\HQ^{#1}(\cX/S)(#2)_{prim}}
\def\HQXp#1#2{\HQ^{#1}(X/S)(#2)_{prim}}
\def\HCcU#1{\HC^{#1}(\cU/S)}
\def\HCcX#1{\HC^{#1}(\cX/S)}
\def\HCcU#1{\HC^{#1}(\cU/S)}
\def\HCcXp#1{\HC^{#1}(\cX/S)_{prim}}
\def\HCXp#1{\HC^{#1}(X/S)_{prim}}
\def\HOcU#1{\HO^{#1}(\cU/S)}
\def\HOcX#1{\HO^{#1}(\cX/S)}
\def\HOcU#1{\HO^{#1}(\cU/S)}
\def\HOcXp#1{\HO^{#1}(\cX/S)_{prim}}
\def\HOXp#1{\HO^{#1}(X/S)_{prim}}
\def\HX#1#2{H^{#1}(X,{#2})}
\def\HXt#1#2#3{H^{#1}(X,{#2}(#3))}
\def\WSan#1{\Omega_{\San}^{#1}}
\def\naban{\nabla^{an}}
\def\HOcUan#1{\HO^{#1}(\cU/S)^{an}}
\def\HOcXan#1{\HO^{#1}(\cX/S)^{an}}
\def\HOXan#1{\HO^{#1}(X/S)^{an}}
\def\HOZan#1{\HO^{#1}(Z/S)^{an}}
\def\HOUan#1{\HO^{#1}(X/S)^{an}}
\def\HOWan#1{\HO^{#1}(W/S)^{an}}
\def\HOYan#1{\HO^{#1}(Y/S)^{an}}
\def\XI{X_I}
\def\Xx{X_x}
\def\Zx{Z_x}
\def\Zjx{Z_{jx}}
\def\Yx{Y_x}
\def\Yjx{Y_{jx}}
\def\Ux{U_x}
\def\Vx{V_x}
\def\xb{\overline{x}}
\def\Uxb{U_{\xb}}

\def\Ut{U_t}
\def\Xt{X_t}
\def\Zt{Z_t}

\def\Zxx{Z_{x}}
\def\Zst{Z}
\def\cZst{\cZ}
\def\Zast{Z^{(\alpha)}}
\def\Wst{W}
\def\Yst{Y}

\def\roo{\kappa_0}
\def\rocXZ{\kappa_{(\cX,\cZ)}}
\def\roX{\kappa_X}
\def\roXZ{\kappa_{(X,Z)}}
\def\roolog{\kappa_0^{log}}
\def\rolog{\kappa^{log}}
\def\rox{\kappa_x}
\def\roxlog{\kappa_x^{log}}

\def\TS{\Theta_S}
\def\TxS{T_x S}
\def\ToS{T_0 S}
\def\UC{U_{\mathbb C}}
\def\XC{X_{\mathbb C}}
\def\SC{S_{\mathbb C}}

\def\TSo{T_0(S)}
\def\WS#1{\Omega_S^{#1}}
\def\WSo#1{\Omega_{S,0}^{#1}}
\def\WX#1{\Omega_X^{#1}}

\def\TX{T_X}
\def\TXS{T_{X/S}}
\def\TXI{T_{\XI}}
\def\TXZx{T_{\Xx}(-\log \Zxx)}
\def\TXx{T_{\Xx}}

\def\WXZ#1{\Omega_X^{#1}(\log \Zst)}
\def\WXZa#1{\Omega_X^{#1}(\log \Zast)}
\def\WcXZ#1{\Omega_{\cX/S}^{#1}(\log \cZst)}
\def\WcX#1{\Omega_{\cX/S}^{#1}}
\def\WcXkZ#1{\Omega_{\cX/k}^{#1}(\log \cZst)}
\def\WcXCZ#1{\Omega_{\cX/\mathbb C}^{#1}(\log \cZst)}
\def\WPnY#1{\Omega_{\mathbb P^n}^{#1}(\log \Yst)}
\def\TXZ{T_X(-\log \Zst)}
\def\TcXZS{T_{\cX/S}(-\log \cZst)}

\def\fXZ#1{\phi_{X,Z}^{#1}}
\def\cXZ{\psi_{(X,Z)}}
\def\cXI{\psi_{\XI}}
\def\cXZx{\psi_{(\Xx,\Zx)}}
\def\cXx{\psi_{\Xx}}
\def\cXo{\psi_{X}}
\def\eXZ{\eta_{(X,Z)}}
\def\ccXS{c_S(\cX)}
\def\ccZIS{c_S(\cZ_I)}
\def\ccXIS{c_S(\cX_I)}
\def\ccXZS{c_S(\cX,\cZ)}
\def\ccXZSo{c_S(\cX,\cZ)_0}
\def\ccXZSx{c_S(\cX,\cZ)_x}
\def\ccXZSt{c_S(\cX,\cZ)_t}
\def\ccXZtS{c_{\wtd{S}}(\wtd{\cX},\wtd{\cZ})}
\def\tS{\wtd{S}}
\def\dpXZ#1{d\rho^{#1}_{X,Z}}
\def\scs{\hbox{ }:\hbox{ }}
\def\onab{\overline{\nabla}}
\def\tonab{@>\onab>>}
\def\OPn#1{\mathcal O_{\bP^n}(#1)}
\def\dlog#1{\frac{d#1}{#1}}
\def\chB#1#2#3{ch_{B,#3}^{#1,#2}}
\def\chD#1#2#3{ch_{D,#3}^{#1,#2}}
\def\chet#1#2#3{ch_{et,#3}^{#1,#2}}
\def\chcont#1#2#3{ch_{cont,#3}^{#1,#2}}
\def\reget#1#2{reg_{et,#1}^{#2}}
\def\regD#1{reg_{D,#1}}
\def\regcont#1{reg_{cont,#1}}
\def\regU#1{reg_U^{#1}}
\def\regUt#1{reg_{U_t}^{#1}}
\def\regetU#1{reg_{et,U}^{#1}}
\def\Ql{\mathbb Q_{\ell}}
\def\Zl{\mathbb Z_{\ell}}
\def\etab{\overline{\eta}}
\def\Qb{\overline{\mathbb Q}}
\def\Pd{\overset{\vee}{\mathbb P}}
\def\kb{\overline{k}}
\def\Ub{\overline{U}}
\def\Sb{\overline{S}}
\def\xb{\overline{x}}
\def\Lb{\overline{L}}
\def\pitop{\pi_1^{top}}
\def\pialg{\pi_1^{alg}}
\def\MB{M_B}
\def\Met{M_{et}}
\def\MBn{M_{B,\nu}}
\def\Metn{M_{et,\nu}}
\def\wU{\omega_U(\sigma)}
\def\wcU{\omega_{\cU/S}(\sigma)}
\def\dcU{\delta_{\cU/S}(\sigma)}
\def\ccU{c_{\cU/S}(\sigma)}
\def\cUx{c_{\Ux}(\sigma)}
\def\cUo{c_{U}(\sigma)}
\def\cUetab{c_{U_{\etab}}(\sigma)}
\def\Dqi{\Delta_{\gamma_i}}
\def\Dq{\Delta_{\gamma}}
\def\nabmo{\overline{\nabla}^{m,0}}
\def\nabmoT{\overline{\nabla}_{0,T}^{m,0}}
\def\nab#1#2{\overline{\nabla}^{#1,#2}}
\def\nabo#1#2{\overline{\nabla}_{0}^{#1,#2}}
\def\naboT#1#2{\overline{\nabla}_{0,T}^{#1,#2}}
\def\otk{\otimes_k}
\def\Pol{P}
\def\h#1#2{h_{#1}(#2)}
\def\San{S_{an}}
\def\Szar{S_{zar}}
\def\SNL{S_{NL}}
\def\SNLq{S_{NL}^q}
\def\dt{\partial_{tame}}
\def\dd{\partial_{div}}
\def\du{\partial_U}
\def\dz{\partial_Z}
\def\uc{\underline{c}}
\def\Tpqc{T_{(p,q)}(\uc)}
\def\Tpqcd{T_{(p',q')}(\uc')}
\def\Tspqc{T^{\sigma}_{(p,q)}(\uc)}
\def\Tspqcd{T^{\sigma}_{(p',q')}(\uc')}
\def\PG#1{\frak{S}_{#1}}

\setcounter{section}{0}
\begin{center}{\bf\arabic{section}. Introduction}\end{center}
\medbreak

In his lectures in [G1], M. Green gives a lucid explanation how fruitful the
infinitesimal method in Hodge theory is in various aspects of algebraic
geometry. A significant idea is to use Koszul cohomology for Hodge-theoretic
computations. The idea originates from Griffiths work [Gri] where the
Poincar\'e residue representation of the cohomology of a hypersurface played
a crucial role in proving the infinitesimal Torelli theorem for hypersurfaces.
Since then many important applications of the idea have been made in different
geometric problems such as the generic Torelli problem and the
Noether-Lefschetz theorem for Hodge cycles and the study of algebraic cycles
(see [G1, Lectures 7 and 8]).
\vskip 3pt

In this paper we apply the method to study an analog of the Noether-Lefschetz
theorem in the context of Beilinson's Hodge conjecture.
Beilinson's Hodge conjecture and its Tate variant concern the regulator
maps for open varieties (cf. [J1, Conjecture 8.5 and 8.6]).
To be more precise we let $U$ be a smooth variety over a field $k$
of characteristic zero.

\begin{description}
\item[(0-1)] (Hodge version)
When $k=\mathbb C$ we have the regulator map from the higher Chow group to
the singular cohomology of $U(\mathbb C)$ (cf. [Bl] and [Sch])
$$ \regU q\scs CH^q(U,q)\otimes\mathbb Q \to
H^q(U,\mathbb Q(q))\cap F^qH^q(U,\mathbb C)$$
where $\mathbb Q(q)=(2\pi\sqrt{-1})^q\mathbb Q$ and $F^*$ denotes the Hodge
filtration of the mixed Hodge structure on the singular cohomology defined by
Deligne [D1].
Taking a smooth compactification $U\subset X$ with $Z=X\setminus U$,
a simple normal crossing divisor on $X$, we have the following formula for
the value of $\regU q$ on decomposable elements in $CH^q(U,q)$;
$$ \regU q(\{g_1,\dots,g_q\})=
\dlog {g_1}\wedge\cdots \wedge \dlog {g_q}\in H^0(X,\WXZ q)=F^qH^q(U,\mathbb C),$$
where $\{g_1,\dots,g_q\}\in CH^q(U,q)$ is the products of
$g_j\in CH^1(U,1)=\Gamma(U,\cO^*_{Zar})$.
\item[(0-2)] (Tate version)
When $k$ is a finitely generated field over $\mathbb Q$ we have the \'etale
regulator map from the higher Chow group to the $\ell$-adic etale cohomology
$$ \regetU q \scs CH^q(U,q)\otimes\Ql \to
H^q_{et}(U_{\kb},\Ql(j))^{Gal(\kb/k)}$$
where $U_{\kb}=U\times_k \kb$ with $\kb$, an algebraic closure of $k$.
\end{description}

\begin{conjecture}\label{conj0-1}
In case $k=\mathbb C$, $\regU q$ is surjective.
\end{conjecture}

\begin{conjecture}\label{conj0-2}
In case $k$ is finitely generated over $\mathbb Q$, $\regetU q$
is surjective.
\end{conjecture}
\rm\medbreak

We call them Beilinson's Hodge and Tate conjectures respectively.
They are analogs of the Hodge and Tate conjectures for algebraic cycles on
smooth projective varieties.
The following are some remarks on the conjectures.
\begin{description}
\item[$(i)$]
The conjectures hold in case $q=1$. They follows from the injectivity
of Abel-Jacobi map (and its $\ell$-adic variant) for divisors on smooth
proper varieties.
\item[$(ii)$]
If $U$ is proper over $k$, we have
\begin{align*}
&H^q(U,\mathbb Q(q))\cap F^qH^q(U,\mathbb C)=0
\quad \hbox{ in case $k=\mathbb C$},\\
&H^q_{et}(U_{\kb},\Ql(q))^{Gal(\kb/k)}=0
\quad \hbox{ in case $k$ is finitely generated over $\mathbb Q$},\\
\end{align*}
by the Hodge symmetry and by the reason of weight respectively.
Thus Beilinson's conjectures are interesting only for non-proper varieties.
\item[$(iii)$]
When $X$ is a projective smooth surface and $U\subset X$ is the complement
of a simple normal crossing divisor $Z\subset X$, then the surjectivity of
$\regU 2$ or $\regetU 2$ has an implication on the injectivity of the regulator
maps for $CH^2(X,1)$, which can be viewed as an analog of the Abel's
theorem for $K_1$ of surfaces. The detail will be explained in \S6.
\item[$(iv)$]
As a naive generalization of Beilinson's Hodge conjecture, one may expect the
surjectivity of the more general regulator map (in case $k=\mathbb C$)
$$ reg_U^{p,q} \scs CH^q(U,2q-p)\otimes\mathbb Q \to
(2\pi\sqrt{-1})^q W_{2q} H^p(U,\mathbb Q)\cap F^qH^p(U,\mathbb C)$$
where $W_*$ denotes the weight filtration.
Jannsen ([J1, 9.11]) has shown that the above map in case $p=1$ and
$q\geq3$ is odd, is not surjective in general by using a theorem of
Mumford [Mu], which implies the Abel-Jacobi for cycles of codimension$\geq 2$
on smooth projective varieties is not injective even modulo torsion.
\end{description}
\medbreak

Before introducing the setup of the main result of this paper, we explain its
background, the Noether-Lefschetz problem for Hodge cycles.
Let $X\subset \mathbb P^n$ be a smooth projective variety over $\mathbb C$.
Recall that the Hodge conjecture predicts that
the space of Hodge cycles in codimension $q$ on $X$:
$$
F^0 H^{2q}(X,\mathbb Q(q)):= H^{2q}(X,\mathbb Q(q)) \cap F^q H^{2q}(X,\mathbb C)$$
is generated by classes of algebraic subvarieties on $X$.
One defines the space of trivial cycles in codimension $q$ on $X$ to be
$$
H^{2q}(X,\mathbb Q(q))_{triv}:=
\Image( H^{2q}(\bP^n,\mathbb Q(q)) \to H^{2q}(X,\bQ(q)))
\subset F^0 H^{2q}(X,\bQ(q)).
$$
It is generated by the class of the section on $X$ of a linear subspace of
codimension $q$ in $\mathbb P^n$.
Now let $S$ be a non-singular quasi-projective variety over $\mathbb C$ and
assume that we are given $\cX \hookrightarrow \mathbb P^n_S$,
an algebraic family over $S$ of smooth projective varieties.
Let $\Xt$ be the fiber of $\cX$ over $t\in S$. Then
the Noether-Lefschetz locus for Hodge cycles in codimension $q$ on $\cX/S$
is defined to be
$$\SNLq=\{t\in S|\; F^0 H^{2q}(\Xt,\mathbb Q(q))\not=
H^{2q}(\Xt,\mathbb Q(q))_{triv}\}.$$
It is the locus of such $t\in S$ that
there are non-trivial Hodge cycles in codimension $q$ on $\Xt$
and hence that the Hodge conjecture is non-trivial for $\Xt$.
One can prove $\SNLq$ is the union of countable number of
(not necessarily proper) closed algebraic subsets of $S$ (cf. [CDK]).
\par

The celebrated theorem of Noether-Lefschetz affirms that in case $\cX/S$
is the universal family of smooth surface of degree$\geq 4$ in $\mathbb P^3$,
$\SNL^1\not=S$, which implies the Picard group of general surface of that kind
is generated by the class of a hyperplane section.
\medbreak

In this paper we propose the following analog of the above problem in the
context of Beilinson's Hodge conjecture.
Assume that we are given schemes over $\mathbb C$
\begin{description}
\item[(0-3)]
$Y=\underset{1\leq j\leq s}{\cup} Y_j \hookrightarrow \mathbb P^n
\hookleftarrow X \hookleftarrow Z=\underset{1\leq j\leq s}{\cup} Z_j,
\quad V=\mathbb P^n\setminus Y,
\quad U=X \setminus Z$
\end{description}
where $X$ is projective smooth, $Y_j\subset \mathbb P^n $ is a smooth
hypersurface,  $Y$ is a simple normal crossing divisor on $\mathbb P^n$,
$Z_j=X \cap Y_j$ intersecting transversally, and $Z$ is a simple normal
crossing divisor on $X$. We are interested in \it the space of
Beilinson-Hodge cycles in degree $q$ on $U$:\rm
$$ F^0 H^q(U,\mathbb Q(q)):= H^q(U,\mathbb Q(q)) \cap F^q H^q(U,\mathbb C)$$
and the regulator map
$\regU q : CH^q(U,q)\otimes\mathbb Q\to F^0 H^q(U,\mathbb Q(q))$.
Since $U$ is affine, $H^q(U,\mathbb Q(q))=0$ for $q>\dim(U)$ by the weak
Lefschetz theorem. Thus we are interested only in case
$1\leq q\leq \dim(U)$. For a positive integer $q\not=n$ we put
$$ H^q(U,\mathbb Q(q))_{triv}=
\Image( H^q(V,\mathbb Q(q)) \to H^q(U,\mathbb Q(q))).$$
We will define a subspace
$ CH^q(U,q)_{triv}\subset CH^q(U,q)\otQ $
such that (cf. Definition\ref{def2-1})
$$ H^q(U,\mathbb Q(q))_{triv}=  \regU q (CH^q(U,q)_{triv})$$
and hence that $H^q(U,\mathbb Q(q))_{triv} \subset \Image(\regU q)
\subset F^0 H^q(U,\mathbb Q(q))$.
We note that in case $\dim(X)>1$ we have (cf. Lemma\ref{lem2-1})
$$ CH^q(U,q)_{triv}=CH^q(U,q)_{dec}$$
where $ CH^q(U,q)_{dec}\subset CH^q(U,q)\otQ $ is the subspace generated
by decomposable elements, namely products of elements in
$CH^1(U,1)=\Gamma(U,\cO_{S_{Zar}}^*)$.
Now assume that we are given schemes over $S$,
a non-singular quasi-projective variety over $\mathbb C$,
\begin{description}
\item[(0-4)]
$\cY=\underset{1\leq j\leq s}{\cup} \cY_j \hookrightarrow \mathbb P^n_S
\hookleftarrow \cX \hookleftarrow \cZ=\underset{1\leq j\leq s}{\cup} \cZ_j,
\quad \cV=\mathbb P^n_S\setminus \cY,
\quad \cU=\cX \setminus \cZ$
\end{description}
whose fibers satisfies the same conditions as (0-3).
Let $\Ut$ be the fiber of $\cU$ over $t\in S$.

\begin{definition}\label{def0-1}
\it For $1\leq q \leq \dim(\Ut)$ we define the Noether-Lefschetz
locus for Beilinson-Hodge cycles in degree $q$ on $\cU/S$ to be
$$ \SNLq=\{x\in S|\; F^0 H^q(\Ut,\mathbb Q(q))\not= H^q(\Ut,\mathbb Q(q))_{triv}\}.$$
\end{definition}
\rm

By definition $\SNLq$ is the locus of such $t\in S$ that
there are non-trivial Beilinson-Hodge cycles in degree $q$ on $\Ut$
and hence that Beilinson's Hodge conjecture is non-trivial for $\Ut$.
We note that a standard method in Hodge theory shows that $\SNLq$ is the union
of countable number of (not necessarily proper) closed analytic subsets of $S$.
\medbreak

A main purpose of this paper is to study the Noether-Lefschetz locus in case
the fibers of $\cX/S$ are smooth complete intersection of
multi-degree $(d_1,\dots,d_r)$ in $\mathbb P^n$.
In this case the Lefschetz theory implies $\SNLq=\emptyset$ unless
$q=m:=n-r=\dim(\Ut)$ (cf. Lemma\ref{lem2-1}).
Thus we are interested only in $\SNL:=\SNL^m$.
In order state the main result we will introduce an invariant $\ccXZS$ in \S3
that measures the ``generality" of the family (0-4), or
how many independent parameters $S$ contains (cf. Definition\ref{def3-1} and
Lemma\ref{lem3-1}).
Assume that the fibers of $\cY_j\subset \mathbb P^n_S$ for $1\leq j\leq s$ is a
hypersurface of degree $e_j$ and put
$$\md=\min\{d_i,e_j|\; 1\leq j\leq s,\; 1\leq i\leq r \}.$$

\begin{theorem}\label{th0-1}
\it For any irreducible component $E\subset \SNL$,
$$\codim_S(E) \geq  \md(n-r-1)+
\underset{1\leq i\leq r}{\sum} d_i  - \ccXZS - n.$$
\end{theorem}
\rm\medbreak

The original proof of the Noether-Lefschetz theorem was based on the study of
monodromy action on cohomology of surfaces (cf. [D2]).
The idea of improving the Noether-Lefschetz theorem by means of the
infinitesimal method in the theory of variations of Hodge structures
was introduced by Carlson, Green, Griffiths and Harris (cf. [CGGH]).
By using the method closer analyses have been made on the Noether-Lefschetz
locus for Hodge cycles on hypersurfaces, particularly on the problem to give a
lower bound of codimensions of its irreducible components and to determine the
components of maximal dimension (cf. [G4], [G5], [V] and [Ot]).
\par

The proof of Theorem\ref{th0-1} follows the same line of arguments as [G4]
by using the variation of mixed Hodge structures arising from cohomology
of the fibers of $\cU/S$. A technical renovation is the results stated in \S1
on generalized Jacobian rings, which give an algebraic description of the
infinitesimal part of the variation of mixed Hodge structures. It is a natural
generalization of the Poincar\'e residue representation of the cohomology of
hypersurfaces in [Gri].
\medbreak

Having a result such as Theorem\ref{th0-1}, a natural quesiton to ask is if the estimate
in the theorem is optimal. In this paper we can give a positive answer to this
quesiton only in case the fibers of $\cX/S$ are plane curves
(cf. Theorem\ref{th4-1}).
In case the fibers of $\cX/S$ are of dimension$>1$, the estimate
seems far from being optimal. Indeed in a forthcoming paper [AS2] we will show
that the optimal estimate of codimensions of the Neother-Lefschetz locus for
Beilinson's Hodge cycles on complements of three hyperplanes in surfaces of
degree $d$ in $\mathbb P^3$ is given explicitly by a quadratic polynomial in
$d$.
\medbreak

Concerning Beilinson's Tate conjecture, we will show the following results.
Let $k$ be a finitely generated field over $\mathbb Q$.
Let $S$ be a quasi-projective variety over $k$ and assume that we are given
schemes over $S$ which satisfy the same condition as (0-4). Assume that the
fibers of $\cX/S$ are smooth complete intersection of multi-degree
$(d_1,\dots,d_r)$ in $\mathbb P^n$. Write $m=n-r$.
\medbreak

\begin{theorem}\label{th0-2}
\it Assume
$\underset{1\leq i\leq r}{\sum} d_i \geq n+1+\ccXZS$.
\begin{description}
\item[(1)]
Let $K=k(S)$ be the function field of $S$. For any a finite generated field
$L$ over $K$ we have
$$ H^m_{et}(U_{\Lb},\Ql(m))^{Gal(\Lb/L)}=
\reget {U_L} m(CH^m(U_L,m)_{triv}\otimes\Ql),$$
where
$ \reget {U_L} m : CH^m(U_L,m)\otimes\Ql \to
H^m_{et}(U_{\Lb},\Ql(m))^{Gal(\Lb/L)}$
is the etale regulator map for $U_L=\cU\times_S \Spec(L)$ (cf. (0-2)).
\item[(2)]
Assume that $k$ is a finite extension of $\mathbb Q$ and
$S(k)\not=\emptyset$.
Let $\pi: S \to \mathbb P_k^N$ be a dominant quasi-finite morphism.
There exist a subset $H\subset\mathbb P_k^N(k)$ such that:
\begin{description}
\item[$(i)$]
For $\forall x\in S$ such that $\pi(x)\in H$ and for any subgroup
$G \subset \Gal(k(\xb)/k(x))$ of finite index we have
$$H^m_{et}(\Uxb,\Ql(m))^G=\reget {\Ux} m(CH^m(\Ux,m)_{triv}\otimes\Ql),$$
where
$ \reget {\Ux} m : CH^m(\Ux,m)\otimes\Ql \to
H^m_{et}(\Uxb,\Ql(m))^{Gal(k(\xb)/k(x))}$
is the etale regulator map for $\Ux$, the fiber of $\cU$ over $x\in S$,
with $\xb$, a geometric point over $x$.
\item[$(ii)$]
Let $\Sigma$ be any finite set of primes of $k$ and let $k_v$ be the
completion of $k$ at $v\in \Sigma$. Then the image of $H$ in
$\prod_{v\in \Sigma} \mathbb P_k^N(k_v)$ is dense.
\end{description}
\end{description}
\end{theorem}
\rm\medbreak

Now we explain how the paper is organized.
In \S1 we state the fundamental results on the generalized Jacobian rings,
the duality theorem and the symmetrizer lemma.
The proof is given in another paper [AS1]. It is based on the basic techniques
to compute Koszul cohomology developed by M. Green ([G2] and [G3]).
In \S2 we define the trivial part of cohomology of complements of
unions of hypersurface sections in smooth projective varieties,
which is necessary to set up the Noether-Lefschetz problem for
Beilinson-Hodge cycles.
In \S3 Hodge theoretic implications of the results in \S1 are stated, which
plays a crucial role in the proof of Theorem\ref{th0-1} given in this section.
The case of plane curves is treated in \S4, where the estimate in
Theorem\ref{th0-1}
is shown to be optimal. Theorem\ref{th0-2} is proven in \S5 by using the results in \S2.
In \S6 we explain an implication of the Beilinson's conjectures on the
injectivity of Chern class maps for $K_1$ of surfaces.
\vskip 20pt

\setcounter{section}{1} \setcounter{subsection}{0}
\begin{center}{\bf\arabic{section}.
Jacobian rings for open complete intersectionn}\end{center}
\medbreak

The purpose of this section is to introduce Jacobian rings for
open compolete intersections and state their fundamental properties.
Throughout the whole paper, we fix integers $r,s \geq 0$ with $r+s\geq 1$,
$n\geq 2$ and $d_1, \cdots, d_r$, $e_1, \cdots, e_s \geq 1$. We put
$$
\d=\sum_{i=1}^{r}d_i, \quad \e=\sum_{j=1}^{s}e_j,\quad
\md=\underset{\underset{1\leq j\leq s}{1\leq i\leq r}}{\min}\{d_i,e_j\},
\quad d_{max}=\underset{1\leq i\leq r}{\max}\{d_i\},
\quad e_{max}=\underset{1\leq j\leq s}{\max}\{e_j\}.
$$
We also fix a field $k$ of characteristic zero.
Let $\Pol=k[X_0,\dots,X_n]$ be the polynomial ring over $k$ in $n+1$
variables.
Denote by $\Pol^l\subset \Pol$ the subspace of the homogeneous
polynomials of degree $l$.
Let $A$ be a polynomial ring over $\Pol$ with indeterminants
$\mu_1, \cdots, {\mu}_r$, $\lambda_1,\cdots, \lambda_s$.
We use the multi-index notation
$$\mu^{\underline{a}}=\mu_1^{a_1}\cdots\mu_r^{a_r},\;
 \lambda^{\underline{b}}=\lambda_1^{b_1}\cdots\lambda_s^{b_s}, \;
\text{for }
\underline{a}=(a_1,\cdots,a_r) \in {\mathbb Z}^{\oplus r}_{\geq 0},\;
\underline{b}=(b_1,\cdots,b_s) \in {\mathbb Z}^{\oplus s}_{\geq 0}.$$
For $q \in {\mathbb Z}$ and $\ell \in {\mathbb Z}$, we write
$$
 A_q(\ell)=
\us{\ba+\bb=q}{\oplus}
\Pol^{\ua\ud+\ub\ue+\ell} \cdot
\mu^{\underline{a}}\lambda^{\underline{b}}\quad
(\ba=\sum_{i=1}^r a_i,\;
\bb=\sum_{j=1}^s b_j,\;
\ua\ud=\sum_{i=1}^r a_id_i,\;
\ub\ue=\sum_{j=1}^s b_je_j)
$$
By convention $A_q(\ell)=0$ if $q<0$.

\begin{definition}\label{def1-1}
For $\ul{F}=(F_1, \cdots, F_r)$, $\ul{G}=(G_1, \cdots, G_s)$ with
$F_i \in \Pol^{d_i}$, $G_j \in \Pol^{e_j}$, we define the {\bf Jacobian ideal}
$J(\ul{F},\ul{G})$ to be the ideal of $A$ generated by
$$\sum_{1\leq i\leq r}\frac{\partial F_i}{\partial X_k}\mu_i+
\sum_{1\leq j\leq s}\frac{\partial G_j}{\partial X_k}\lambda_j,
\quad F_i,\quad G_j\lambda_j \quad
(1 \leq i \leq r,\; 1 \leq j \leq s,\; 0 \leq k \leq n).$$
The quotient ring $B=B(\ul{F},\ul{G})=A/J(\ul{F},\ul{G})$ is called
the {\bf Jacobian ring of} ($\ul{F},\ul{G}$). We denote
$$
B_q(\ell) =B_q(\ell)(\ul{F},\ul{G})=A_q(\ell)/A_q(\ell) \cap J(\ul{F},\ul{G}).
$$
\end{definition}
\medbreak

\begin{definition}\label{def1-2}
Suppose $n\geq r+1$.
Let $\mathbb{P}^n=\text{\rm Proj }\Pol$ be the projective space over $k$.
Let $X\subset \P^n$ be defined by $F_1=\cdots=F_r=0$ and let
$Z_j\subset X$ be defined by $G_j=F_1=\cdots=F_r=0$ for $1 \leq j \leq s$.
We also call $B(\ul{F},\ul{G})$ the Jacobian ring of the pair
$(X,\Zst={\cup}_{1\leq j\leq s} Z_j)$ and denote
$B(\ul{F},\ul{G})=B(X,\Zst)$ and $J(\ul{F},\ul{G})=J(X,\Zst)$.
\end{definition}
\medbreak

In what follows we fix $\ul{F}$ and $\ul{G}$ as Definition\ref{def1-1} and
assume the condition
\begin{description}
\item[(1-1)]
$F_i=0$ ($1\leq i\leq r$) and
$G_j=0$ ($1\leq j\leq s$) intersect transversally in $\mathbb{P}^n$.
\end{description}
We mension three main theorems.
The first main theorem concerns with the geometric meaning of Jacobian rings.

\begin{theorem}\label{thI}
\it
Suppose $n\geq r+1$. Let $X$ and $Z$ be as Definition\ref{def1-2}.
\begin{description}
\item[(1)]
For intergers $0 \leq q\leq n-r$ and $\ell\geq 0$ there is a natural
isomorphism
$$\fXZ q \scs B_q(\d+\e-n-1+\ell) \isom H^{q}(X, \WXZ {n-r-q}(\ell))_{prim}.$$
Here $\WXZ p$ is the sheaf of algebraic differential $q$-forms on $X$
with logarithmic poles along $\Zst$ and `$prim$' means the primitive part:
$$ H^{q}(X, \WXZ {p}(\ell))_{prim}=
\left.\left\{\begin{gathered}
   \Coker(H^q(\bP^n,\Omega_{\bP^n}^q) \to H^{q}(X, \WX {q}))\\
   H^{q}(X, \WXZ p(\ell))\\
\end{gathered}\right.
\begin{aligned}
  &\text{ if $q=p$ and $s=\ell=0$,}\\
  &\text{ otherwise.}\\
\end{aligned}\right.
$$
\item[(2)]
There is a natural map
$$\cXZ : B_1(0) \longrightarrow H^1(X, \TXZ)_{alg}\subset H^1(X, \TXZ) $$
which is an isomorphism if $\dim(X)\geq 2$.
Here $\TXZ$ is the $\cO_X$-dual of $\WXZ 1$ and the group on the right hand
side is defined in Definition\ref{def1-3} below. The following map
$$
H^1(X, \TXZ) \otimes H^{q}(X, \WXZ p)
\longrightarrow H^{q+1}(X, \WXZ {p-1}).
$$
induced by the cup-product and the contraction
$\TXZ\ot\WXZ p\to\WXZ {p-1}$
is compatible through $\cXZ$ with the ring multiplication up to scalar.
\end{description}
\end{theorem}
\rm\medbreak

Roughly speaking, the generalized Jacobian rings describe the infinitesimal
part of the Hodge structures of open
variety $X \setminus Z$,
and the cup-product with Kodaira-Spencer class
coincides with the ring multiplication up to non-zero scalar.
This result was originally invented by P. Griffiths in case of hypersurfaces
and generalized to complete intersections by Konno [K].
Our result is a further generalization.

\begin{definition}\label{def1-3}
Let the assumption be as in Theorem\ref{thI}.
We define $H^1(X, \TXZ)_{alg}$ to be
the kernel of the composite map
$$H^1(X,\TXZ)\to H^1(X,\TX) \to H^2(X,\cO_X),$$
where the second map is induced by the cup product with the class
$c_1(\cO_{X}(1))\in H^1(X,\WX 1)$ and the contraction
$\TX\otimes \WX 1 \to \cO_X$.
It can be seen that
$$ \dim_k(H^1(X, \TXZ)/H^1(X, \TXZ)_{alg})=
\left.\left\{\begin{gathered}
   1 \\
   0 \\
\end{gathered}\right.
\begin{aligned}
  &\text{ if $X$ is a $K3$ surface,}\\
  &\text{ otherwise.}\\
\end{aligned}\right.
$$
\end{definition}
\rm\medbreak

The second main theorem is the duality theorem for the generalized
Jacobian rings.

\begin{theorem}\label{thII}
\it Let the notation be as as above.
\begin{description}
\item[(1)]
There is a natural map (called the trace map)
$$\tau\;:\; B_{n-r}(2(\d-n-1)+\e) \to k.$$
Let
$$ \h p \ell\;:\; B_p(\d-n-1+\l) \ra B_{n-r-p}(\d+\e-n-1-\l)^*$$
be the map induced by the following pairing induced by
the multiplication
$$B_{p}(\d-n-1+\l) \otimes B_{n-r-p}(\d+\e-n-1-\l)\to  B_{n-r}(2(\d-n-1)+\e)
\rmapo{\tau} k.$$
When $r>n$ we define $ \h p \ell$ to be the zero map by convention.
\item[(2)]
The map $\h p \ell$ is an isomorphism in either of the following cases.
\begin{description}
\item[$(i)$]  $s\geq 1$ and $p<n-r$ and $\ell<e_{max}$.
\item[$(ii)$]  $s\geq 1$ and $0\leq \ell\leq e_{max}$ and $r+s\leq n$.
\item[$(iii)$] $s=\ell=0$ and either $n-r\geq 1$ or $n-r=p=0$.
\end{description}
\item[(3)]
The map $\h {n-r} \ell$ is injective if $s\geq 1$ and $\ell<e_{max}$.
\end{description}
\end{theorem}
\rm\medbreak

We have the following auxiliary result on the duality.

\begin{theorem}\label{thII'}
\it
Assume $n-r\geq 1$ and consider the composite map
$$
\eXZ \scs H^{0}(X, \WXZ {n-r})\rmapo{(\fXZ 0)^{-1}}
B_0(\d+\e-n-1) \rmapo{{\h {n-r} 0}^*}  B_{n-r}(\d-n-1)^*
$$
where the second map is the dual of $\h {n-r} 0$.
Then $\eXZ$ is surjective and we have (cf. Definition\ref{def1-4} below)
$$ \Ker(\eXZ)=\wedge_X^{n-r}(G_1,\dots,G_s).$$
\end{theorem}
\rm\medbreak

\begin{definition}\label{def1-4}
Let $G_1,\dots,G_s$ be as in Definition\ref{def1-1} and let $Y_j\subset \P^n$
be the smooth hypersurface defined by $G_j=0$.
Let $X\subset \P^n$ be a smooth projective variety such that
$Y_j$ ($1\leq j\leq s$) and $X$ intersect transversally.
Put $Z_j=X\cap Y_j$.
Take an integer $q$ with $0\leq q \leq s-1$.
For integers $1\leq j_1<\cdots< j_{q+1}\leq s$, let
$$ \omega_X(j_1,\dots,j_{q+1})\in H^{0}(X, \WXZ {q})
\quad (\Zst=\underset{1\leq j\leq s}{\Sigma} Z_j)$$
be the restriction of
$$ \sum_{\nu=1}^{q+1}(-1)^{\nu-1} e_{j_{\nu}}
\frac{dG_{j_1}}{G_{j_1}}\wedge\cdots\wedge\widehat
{\frac{dG_{j_\nu}}{G_{j_\nu}}}
\wedge\cdots\wedge{\frac{dG_{j_{q+1}}}{G_{j_{q+1}}}} \;\in
H^0(\P^n,\WPnY q)$$
where
$\Yst=\underset{1\leq j\leq s}{\Sigma} Y_j\subset \P^n$.
We let
$$ \wedge_X^q(G_1,\dots,G_s)\subset H^{0}(X, \WXZ {q})$$
be the subspace generated by $\omega_X(j_1,\dots,j_{q+1})$.
For $1\leq j_1<\cdots< j_{q}\leq s-1$ we have
$$  e_1^{q-1}\cdot \omega_X(1,j_1,\dots,j_{q})=
\frac{dg_{j_1}}{g_{j_1}}\wedge\cdots\wedge
{\frac{dg_{j_{q}}}{g_{j_{q}}}}
\qwith g_j=(G_j^{e_1}/G_1^{e_j})_{|X}\in \Gamma(U,\cO_{U}^*)
\quad (U=X\setminus Z)$$
and $\omega_X(1,j_1,\dots,j_{q})$ with $1<j_1<\cdots< j_{q}\leq s$
form a basis of $\wedge_X^q(G_1,\dots,G_s)$.
\end{definition}
\rm\medbreak

Our last main theorem is the generalization of Donagi's symmetrizer lemma [Do]
(see also [DG], [Na] and [N]) to the case of open complete intersections at
higher degrees.

\begin{theorem}\label{thIII}
\it Assume $s\geq 1$.
Let $V \subset B_1(0)$ is a subspace of codimension $c\geq 0$.
Then the Koszul complex
$$
B_p(\l) \otimes \os{q+1}{\wedge}V \rightarrow
B_{p+1}(\l) \otimes \os{q}{\wedge}V \rightarrow
B_{p+2}(\l) \otimes \os{q-1}{\wedge}V
$$
is exact if one of the following conditions is  satisfied.
\begin{description}
\item[$(i)$]
$p\geq 0$, $q=0$ and $\md p+\ell\geq c$.
\item[$(ii)$]
$p\geq 0$, $q=1$ and $\md p+\ell\geq 1+c$ and $\md(p+1)+\ell\geq d_{max}+c$.
\item[$(iii)$]
$p\geq 0$, $\md(r+p)+\ell\geq \d+q+c$, $\d+e_{max}-n-1>\ell\geq \d-n-1$
and either $r+s\leq n+2$ or $p\leq n-r-[q/2]$, where
$[*]$ denotes the Gaussian symbol.
\end{description}
\end{theorem}
\rm\medbreak
\vskip 20pt

\setcounter{section}{2}\setcounter{subsection}{0}
\begin{center}{\bf\arabic{section}.
Trivial part of cohomology of complements of hypersurfaces}\end{center}
\medbreak

In this section we introduce some notions necessary to set up the
Noether-Lefschetz problem for Beilinson-Hodge cycles on the complement of
the union of hypersurface sections in a smooth projective variety.
We fix the base field $k$ which is either $\mathbb C$ or finitely generated
over $\mathbb Q$. Assume that we are given schemes over $k$
\begin{description}
\item[(2-1)]
$Y=\underset{1\leq j\leq s}{\cup} Y_j \hookrightarrow \mathbb P^n
\hookleftarrow X \hookleftarrow Z=\underset{1\leq j\leq s}{\cup} Z_j,
\quad V=\mathbb P^n\setminus Y,
\quad U=X \setminus Z$
\end{description}
where $X$ is projective smooth, $Y_j\subset \mathbb P^n $ is a smooth
hypersurface of degree $e_j$, $Y$ is a simple normal crossing divisor on
$\mathbb P^n$, $Z_j=X \cap Y_j$ intersecting transversally, and $Z$ is a
simple normal crossing divisor on $X$. In what follows $CH^i(*,j)$ denotes
the Bloch's higher Chow group (cf. [Bl]).

\begin{definition}\label{def2-1}
\it Assume $s\geq 2$.
\begin{description}
\item[(1)]
Let $G_j\in H^0(\mathbb P^n,\cO_{\mathbb P^n}(e_j))$ be a non-zero
element defining $Y_j\subset \mathbb P^n $.
For $1\leq j\leq s-1$ put
$$g_j=(G_j^{e_s}/G_s^{e_j})_{|U}\in \Gamma(U,\cO_{U}^*)=CH^1(U,1).$$
Let
$CH^1(U,1)_{triv}\subset CH^1(U,1)\otimes \mathbb Q$
be the subspace generated by $k^*$ and $g_j$ with $1\leq j\leq s-1$.
\item[(2)]
For $q\geq 1$ let
$CH^q(U,q)_{dec}\subset CH^q(U,q)\otimes \mathbb Q$ be the subspace generated
by the products of elements in $CH^1(U,1)$ and let
$$  CH^q(U,q)_{triv}\subset CH^q(U,q)_{dec} $$
to be the subspace generated by the products of elements in $CH^1(U,1)_{triv}$.
\item[(3)]
For $q\geq 1$ let
$  H^0(X,\WXZ q)_{triv}\subset H^0(X,\WXZ q) $
be the subspace generated by the wedge products of
$dg_j/g_j\in H^0(X,\WXZ 1)$ with $1\leq j\leq s-1$.
\item[(4)]
Let
$ H^p(U,q)$ denotes the singular cohomology
$H^p(U,\mathbb Q(q))$ of $U(\mathbb C)$ in case $k=\mathbb C$
(resp. the \'etale cohomology $H^p_{et}(\Ub,\Ql(q))$
in case $k$ is finitely generated over $\mathbb Q$, where
$\Ub=U\times_k \kb$ with $\kb$, an algebraic closure of $k$).
For $q\geq 1$ let
$$ H^q(U,q)_{triv}\subset H^q(U,q)$$
be the subsapce generated by the image of $CH^q(U,q)_{triv}$
under the regulator map introduced in (0-1) and (0-2).
\end{description}
In case $s=0,1$ we put $M_{triv}=0$ by convention for
$M=CH^q(U,q)$, $H^0(X,\WXZ q)$, $H^q(U,q)$.
\end{definition}
\rm\medbreak

\begin{remark}\label{rem2-1}
\it In case $k=\mathbb C$ we have
$$ H^q(U,\mathbb Q(q))_{triv}\otimes_{\mathbb Q}\mathbb C=
H^0(X,\WXZ q)_{triv}$$
under the identification
$F^q H^q(U,\mathbb C)=H^0(X,\WXZ q)$.
\end{remark}
\rm\medbreak

\begin{lemma}\label{lem2-1}
\it Assume $s\geq 1$ and $q\geq 1$.
\begin{description}
\item[(1)]
If $\dim(X)>1$, $CH^q(U,q)_{dec}=CH^q(U,q)_{triv}$.
\item[(2)]
If $q\not=n$,
$ H^q(U,q)_{triv}=\Image(H^q(V,q) \to H^q(U,q)).$
\item[(3)]
If $X\subset \mathbb P^n$ is a smooth complete intersection,
$ H^q(U,q)=H^q(U,q)_{triv}$ for $q\not=\dim(X)$.
\end{description}
\end{lemma}

\rm
\begin{proof}
To show Lemma\ref{lem2-1}(1) it suffices to show
$CH^1(U,1)_{triv}=CH^1(U,1)\otimes \mathbb Q$,
which follows at once from the exact sequence
$$ 0\to k^* \to CH^1(U,1) \to \underset{1\leq j\leq s}{\bigoplus}
\mathbb Z \to CH^1(X)$$
where the second map is induced by the residue maps along $Z_j$ and
the last by the class of $Z_j$ in $CH^1(X)$.

We show Lemma\ref{lem2-1}(3).
Since $U$ is affine by the assumption $s\geq 1$,
$H^q(U,q)=0$ for $q>\dim(X)$ by the weak Lefschetz theorem.
Thus we may assume that $1\leq q <\dim(X)$ so that $\dim(X)>1$.
If $s=1$, $Z=Z_1$ is a smooth complete intersection of dimension$\geq 1$
and we have the exact sequence
$$ H^{q-2}(Z,q-1) \rmapo{\alpha} H^q(X,q) \to H^q(U,q) \to
 H^{q-1}(Z,q-1) \rmapo{\beta} H^{q+1}(X,q)$$
where $\alpha$ and $\beta$ are the Gysin maps. The assumption $q\not=\dim(X)$
implies $q-1\not=\dim(Z)$ and the Lefschetz theory implies that $\alpha$
is surjective and $\beta$ is injective so that $H^q(U,q)=0$. If $q=1$ the
desired assertion follows from the exact sequence
$$ 0\to H^1(U,1) \to \underset{1\leq j\leq s}{\bigoplus} H^0(Z_j,0)
\to H^2(X,1),$$
where the last map is the Gysin map and the injectivity of the first map
follows from the vanishing of $H^1(X,1)$ by the assumption on $X$.
Now assuming $\dim(X)>1$, $s>1$ and $q>1$, we proceed by the double induction
on $\dim(X)$ and $s$. Put
$$U'=X\setminus \underset{j=2}{\overset{s}{\cup}} Z_j,\quad
W=U'-U=Z_1\setminus (Z_1 \cap
\underset{j=2}{\overset{s}{\cup}} Z_j).$$
We have the short exact sequence
$$ H^q(U',q)\to H^q(U,q) \to H^{q-1}(W,q-1).$$
By induction
$H^q(U',q)=H^q(U',q)_{triv}$ and $H^{q-1}(W,q-1)=H^{q-1}(W,q-1)_{triv}$
while it is easy to check that the residue map
$H^q(U,q)_{triv} \to H^{q-1}(W,q-1)_{triv}$ is surjective.
This proves the desired assertion for $H^q(U,q)$.

To show Lemma\ref{lem2-1}(2) let $H^q(V,q)_{triv}\subset H^q(V,q)$ be defined
as $H^q(U,q)_{triv}$ by taking $X=\mathbb P^n$.
By definition we have
$ H^q(U,q)_{triv}=\Image(H^q(V,q)_{triv} \to H^q(U,q)).$
By Lemma\ref{lem2-1}(3) $H^q(V,q)_{triv}=H^q(V,q)$ if $q\not=n$.
This proves the desired assertion.
\end{proof}

\setcounter{section}{3}\setcounter{subsection}{0}
\begin{center}{\bf\arabic{section}.
Hodge theoretic implication of generalized Jacobian rings}\end{center}
\medbreak

In this section we prove Theorem\ref{th0-1} by using the results in \S1.
Let the assumption be as in \S2.
We fix a non-singular affine algebraic
variety $S$ over $k$ and the following schemes over $S$
\begin{description}
\item[(3-1)]
$\cY=\underset{1\leq j\leq s}{\cup} \cY_j \hookrightarrow \mathbb P^n_S
\hookleftarrow \cX \hookleftarrow \cZ=\underset{1\leq j\leq s}{\cup} \cZ_j,
\quad \cV=\mathbb P^n_S\setminus \cY,
\quad \cU=\cX \setminus \cZ$
\end{description}
whose fibers are as (2-1). We assume that the fibers of $\cX/S$ are smooth
complete intersections of multi-degree $(d_1,\dots,d_r)$ and those of
$\cY_j\subset \mathbb P^n_S$ for $1\leq j\leq s$ are hypersurfaces of degree
$e_j$. For integers $p,q$ we introduce the following sheaf on $\Szar$
$$\HcU p q=R^q f_* \WcXZ p,$$
where $f:\cX\to S$ is the natural morphism and
$\WcXZ p=\os{p}{\wedge}\WcXZ 1$ with $\WcXZ 1$, the sheaf of relative
differentials on $\cX$ over $S$ with logarithmic poles along $\cZ$.
In case $s\geq 1$ the Lefschetz theory implies $\HcU p q=0$
if $p+q\not=n-r$.
In case $s=0$ it implies $\HcXp p q=0$
if $p+q\not=n-r$ where ``$prim$" denotes the primitive part
(cf. Theorem\ref{thI}(1)).
The results in \S1 implies that under an appropriate
numerical condition on $d_i$ and $e_j$ we can control the cohomology of the
following Koszul complex
$$ \WS {q-1}\ot \HcU {a+2}{b-2} \rmapo{\onab}
 \WS {q}\ot \HcU {a+1}{b-1} \rmapo{\onab}
 \WS {q+1}\ot \HcU {a}{b}.$$
Here $\onab$ is induced by the Kodaira-Spencer map
\begin{description}
\item[(3-1)]
$ \rocXZ\scs \TS \to R^1f_* \TcXZS,$
\end{description}
with $\TS=\cHom_{\cO_S}(\WS 1,\cO_S)$ and
$\TcXZS=\cHom_{\cO_{\cX}}(\WcXZ 1,\cO_{\cX})$, and the map
$$ R^1f_* \TcXZS \otimes R^{b-1} f_* \WcXZ {a+1} \to R^b f_* \WcXZ a$$
induced by the cup product and
$\TcXZS \otimes \WcXZ {a+1} \to \WcXZ a$, the contraction.
For the application to the Beilinson's conjectures the kernel of
the following map plays a crucial role
$$  \nab p q \scs \HcU p q \to \WS {1}\otO \HcU {p-1}{q+1}.$$
In case $s=0$ we let $\nab p q$ denote the primitive part of the above map:
$$ \HcXp p q \rmapo{\onab} \WS {1}\otO \HcXp {p-1}{q+1}.$$
Now we fix a $k$-rational point $0\in S$ and let
$U\subset X \supset Z$ denote the fibers of $\cU\subset \cX \supset \cZ$.
Taking the fibers of the above maps we get
\begin{align*}
& \nabo p q \scs \HU p q \to \WSo {1}\otimes \HU {p-1}{q+1}\quad (s\geq 1),\\
&\nabo p q \scs \HXp p q \to \WSo {1}\otimes \HXp {p-1}{q+1}\quad (s=0)\\
\end{align*}
where $\HU p q= H^q(X,\WXZ p)$ and $\HXp p q=H^q(X,\WX p)_{prim}$.
Let $\TSo$ be the tangent space of $S$ at $0$ and fix a linear subspace
$T\subset \TSo$. Via the canonical isomorphism
$ \WSo 1 \cong \Hom_k(\TSo,k)$, the above maps induce
\begin{align*}
& \naboT p q \scs \HU p q \to \Hom_k(T,\HU {p-1}{q+1})\quad (s\geq 1),\\
& \naboT p q \scs \HXp p q \to \Hom_k(T,\HXp {p-1}{q+1})\quad (s=0)\\
\end{align*}
Now the key result is the following. See Definition\ref{def3-1} below for the
definition of $\ccXZSo$. Put
$$\d=\underset{1\leq i\leq r}{\sum} d_i,\quad
\md=\min\{d_i,e_j|\; 1\leq j\leq s,\; 1\leq i\leq r \}.$$

\begin{theorem}\label{th3-1}
\it Assume $p+q=m:=n-r\geq 1$. Let $c=\codim_{\TSo}(T)$.
\begin{description}
\item[(1)]
Assuming $1\leq p\leq m-1$ and
$\md(p-1)+\d\geq n+1+\ccXZSo+c$, $\Ker(\naboT pq)=0.$
\item[(2)]
Assuming $\md(n-r-1)+\d\geq n+1+\ccXZSo+c$,
$\Ker(\nabmoT)=H^0(X,\WXZ m)_{triv}$ (cf. Definition\ref{def2-1}).
\end{description}
\end{theorem}
\rm\medbreak

\begin{definition}\label{def3-1}
\it We define
$$\ccXZSo=\dim_{k}(\Image(\cXZ)/\Image(\cXZ)\cap \Image(\roXZ))\}.$$
where
$\roXZ \scs \TSo \to H^1(X,\TXZ)$
(resp. $\cXZ: B_1(0) \to H^1(X,\TXZ)$)
is the fiber at $0$ of the map (3-2)
(resp. the map in Theorem\ref{thI}(2) for $(X,Z)$). We also define
$$\ccXZS=\underset{\; t\in S}{\min}
\{\ccXZSt\}$$
where $\ccXZSt$ is defined by the same way as $\ccXZSo$ for the fibers of
$\cU\subset\cX\supset\cZ$ over $t\in S$.
\end{definition}
\rm\medbreak

\begin{remark}\label{rem3-1}
If $n-r\geq 2$ and $X$ is not a $K3$ surface,
$\cXZ$ is surjective so that
$\ccXZSo=\dim_k\Coker(\roXZ).$
\end{remark}
\medbreak

Let $P=\mathbb C[z_0,\dots,z_n]$ be the homogeneous coordinate ring of $\mathbb P^n$
and let $P^d\subset P$ be the space of the homogeneous polynomials of degree
$d$. The dual projective space
$$\Pd(P^d)=\mathbb P_k^{N(n,d)}\quad (N(n,d)=\binom{n+d}{d}-1)$$
parametrizes hypersurfaces
$Y\subset \mathbb P^n$ of degree $d$ defined over $k$. Let
$$ B\subset \underset{1\leq i\leq r}{\prod} \mathbb P_k^{N(n,d_i)} \times
\underset{1\leq j\leq s}{\prod} \mathbb P_k^{N(n,e_j)}$$
be the Zariski open subset parametrizing such
$((X_i)_{1\leq i \leq r},(Y_j)_{1\leq j \leq s})$
that $X_1+\cdots + X_r+ Y_1+\cdots+Y_s$ is a simple normal crossing divisor on
$\mathbb P^n_k$. We consider the family
$$\cXB \hookleftarrow \cZB=\underset{1\leq j\leq s}{\cup}
\cZBj \quad\text{ over } B
$$
whose fibers are
$ X \hookleftarrow Z=\underset{1\leq j\leq s}{\cup}Z_j $
with $X=X_1\cap\cdots \cap X_r$ and $Z_j=X\cap Y_j$.
We omit the proof of the following.

\begin{lemma}\label{lem3-1}
\it
Let $E\subset B$ be a non-singular locally closed subvariety of codimension
$c\geq 0$ and let $S\to E$ be a dominant map. Assume that the family
$(\cX,\cZ)/S$ is the pullback of $(\cXB,\cZB)/B$ via $S\to B$.
Then we have $\ccXZS\leq c$.
\end{lemma}
\rm\medbreak

Now we prove Theorem\ref{th3-1}.
We only show the second assertion and leave the first to the readers.
The fact $H^0(X,\WXZ m)_{triv}\subset \Ker(\nabmoT)$ follows from the fact that
it lies in the image of
$$H^0(\cX,\WcXkZ m)\to \Gamma(S,f_*\WcXZ m) \to H^0(X,\WXZ m),$$
where $\WcXkZ \cdot$ is the sheaf of differential forms of $\cX$ over $k$
with logarithmic poles along $\cZ$. Thus we are reduced to show the
injectivity of the induced map
$$H^0(X,\WXZ m)/H^0(X,\WXZ m)_{triv} \to \Hom_k(T,H^1(X,\WXZ {m-1})).$$
By Theorem\ref{thII} and Theorem\ref{thII'}
in \S1 this is reduced to show the surjectivity of
$$ W\otimes B_{n-r-1}(\d-n-1) \to B_{n-r}(\d-n-1)$$
where $B_1(0)\supset W:=\cXZ^{-1}(\roXZ(T))$.
By definition $W$ is of codimension$\leq c+\ccXZSo$ in $B_1(0)$.
Hence the desired assertion follows from Theorem\ref{thIII}$(i)$ in \S1.
\medbreak

Having Theorem\ref{th3-1}, Theorem\ref{th0-1} is proven by a standard method in Hodge theory
(cf. [G1, p.75]), which we recall in what follows. Let the assumption be as in
Theorem\ref{th0-1}. Fix $0\in E\subset \SNL$ and let
$U\subset X \supset Z$ be as before. Choosing a simply connected neighbourhood
$\Delta$ of $0$ in $S(\mathbb C)$ we have the natural identification for
$\forall t\in \Delta$ via flat translation with respect to the Gauss-Maninn
connection:
$$ H^m(U,\mathbb Q(m)) \isom H^m(\Ut,\mathbb Q(m));\; \lambda \to \lambda_t$$
where $\Ut$ is the fiber of $\cU$ over $t\in S$.
By definition of $\SNL$ there exists $\lambda\in F^0H^m(U,\mathbb Q(m))$
not contained in $H^m(U,\mathbb Q(m))_{triv}$
such that $\lambda_t\in F^m H^m(\Ut,\mathbb C)$ for $\forall t\in E\cap \Delta$.
By noting Remark\ref{rem2-1} it implies that
$\lambda\in \Ker(\nabmoT)\not=H^0(X,\WXZ m)_{triv}$,
where $T\subset \TSo$ is the tangent space of $E$. Hence we have
$$ \codim_S(E)\geq \codim_{\TSo}(T)\geq \md(n-r-1)+
\underset{1\leq i\leq r}{\sum} d_i -n-\ccXZS$$
where the second inequality follows from Theorem\ref{th3-1}.
This completes the proof of Theorem\ref{th0-1}.

\vskip 20pt

\setcounter{section}{4}\setcounter{subsection}{0}
\begin{center}{\bf\arabic{section}.
Case of plane curves}\end{center}
\medbreak

One may naturally asks if the estimate in Theorem\ref{th0-1} is optimal.
In this section we answer the question in case $n=2$ and $r=1$, namely
in case the fibers of $\cX/S$ are plane curves. We work over $\mathbb C$.
Let $P=\mathbb C[z_0,z_1,z_2]$ be the homogeneous coordinate ring of $\mathbb P^2$
and let $P^d\subset P$ be the space of the homogeneous polynomials of degree
$d$.
Let
$$ S \subset  P^d/\mathbb C^* \times
\underset{1\leq j\leq s}{\prod} P^{e_j}/\mathbb C^*
= \mathbb P^{N(2,d)} \times
\underset{1\leq j\leq s}{\prod} \mathbb P^{N(2,e_j)}
\quad (N(2,u)=\binom{u+2}{2}-1)$$
be the moduli space of such sets of plane curves $(X,(Y_j)_{1\leq j\leq s})$
that $X \cup Y_1 \cup \cdots \cup Y_s$ is a simple normal crossing divisor on
$\mathbb P^2$. Let
$\cX \hookrightarrow \mathbb P^n_S \hookleftarrow \cY_j$
be the universal families over $S$ and put
$\cZ=\cX\cap \underset{1\leq j\leq s}{\cup} \cY_j$ and
$\cU=\cX\setminus \cZ$. Let $\Ut\subset\Xt\supset\Zt$ denote the fibers of
$\cU\subset\cX\supset\cZ$ over $t\in S$ .
We are interested in the Noether-Lefschetz locus for Beilinson-Hodge cycles
on $\Ut$;
$$ \SNL=\{t\in S|\; F^0H^1(\Ut,\mathbb Q(1))\not=H^1(\Ut,\mathbb Q(1))\}.$$
Noting that $\regUt 1$ is known to be surjective and that
$\Ker(\regUt 1)=\mathbb C^*\otimes\mathbb Q \subset CH^1(\Ut,1)\otimes\mathbb Q,$
we have
$$ \SNL=\{t\in S|\; CH^1(\Ut,1)\otimes\mathbb Q \not=CH^1(\Ut,1)_{triv}\}.$$

\begin{theorem}\label{th4-1}
\it Let $E\subset \SNL$ be an irreducible component. Then
$\codim_S(E)\geq d-2$. Moreover there exists $E$ for which the equality holds.
\end{theorem}

\rm\begin{proof}
The inequality follows from Theorem\ref{th0-1} together with the fact $\ccXZS=0$,
which follows from Lemma\ref{lem3-1}.
As for the second assertion we consider the subset $\Sigma\subset S$ of those
$t\in S$ that satisfy the condition: there exists the unique point $x_t\in \Zt$
such that $\Xt\cap H_t =\{x_t\}$ set-theoretically, where $H_t$ is the tangent
line of $\Xt$ at $x_t$. It is easy to check that $\Sigma$ is locally closed
algebraic subset of $S$. Thus it suffices to show the following.
\end{proof}

\begin{lemma}\label{lem4-1}
\it $\Sigma\subset \SNL$, $\Sigma$ is equidimensional,
and $codim_S(\Sigma)=d-2$.
\end{lemma}

\rm\begin{proof}
Fix any $t\in \Sigma$ and let
$Y_j$ and $U$ denote the fibers of $\cY_j$ and $\cU$ over $t$ for simplicity.
Choose $L\in P^1$ defining $H_t\subset \mathbb P^2$ and $G_j\in P^{e_j}$ defining
$Y_j$ which contains $x_t$. Then it is easy to see that
$(G_j/L^{e_j})_{|U}\in CH^1(U,1)$ while it is not contained in
$CH^1(U,1)_{triv}$. This proves the first assertion of Lemma\ref{lem4-1}.
In order to show the remaining assertions, let
$M$ be the space of pairs $(x,H)$ of a point $x\in \mathbb P^2$ and
a line $H\subset \mathbb P^2$ passing through $x$. Define the morphism
$ \pi\scs \Sigma \to M$ by $\pi(t)=(x_t,H_t)$. Considering the action of
$\PGL_3(\mathbb C)$, the group of linear transformations on $\mathbb P^2$,
one easily see that $\pi$ is surjective.
Thus it suffices to show the following.
\end{proof}

\begin{lemma}\label{lem4-2}
\it Let $\Sigma_0 \subset \Sigma$ be the fiber of $\pi$ over
$(0,H_0)$ with $0:=(0:0:1)$ and $H_0\subset \mathbb P^n$ defined by $z_0\in P^1$.
Then $\Sigma_0$ is equidimensional and $codim_S(\Sigma)=d-2+\dim(M)=d+1$.
\end{lemma}

\rm\begin{proof}
Let $T\subset S':=\underset{1\leq j\leq s}{\prod} \mathbb P_k^{N(2,e_j)}$
be the closed subset parametrizing such $(Y_j)_{1\leq j\leq s}$ that
$0\in Y_j$ for some $j$. Clearly $T$ is equidimensional and of
$\codim_{S'}(T)=1$.
By definition one sees that $\Xt$ for $t\in E_0$ is defined by
an equation of the form $ z_0 A +z_1^d $ with $A\in P^{d-1}$.
This implies that $E_0$ is identified with a non empty open subset of
$P^{d-1}\times T$, where $P^{d-1}$ is veiwed as an affine space over $\mathbb C$.
This shows that $E_0$ is equidimensional and that
$\codim_S(E_0)=\binom{d+2}{2}-1-\binom{d+1}{2}+\codim_{S'}(T)=d+1$ as desired.
\end{proof}
\medbreak

\begin{remark}\label{rem4-1}
It is interesting to ask if the components $E$ of $\Sigma$
are the only ones satisfying $\codim_S(E)=d-2$.
It is possible to give a positive answer to this question in case
$s=2$ and $e_1=e_2=1$ by using the method in [Gr5] and [V].
The general case remains open.
\end{remark}

\vskip 20pt

\setcounter{section}{5}\setcounter{subsection}{0}
\begin{center}{\bf\arabic{section}.
Beilinson's Tate conjecture}\end{center}
\medbreak

In this section we show Theorem\ref{th0-2}. Let the assumption be as in the beginning of
\S3. Write $m=n-r$ and $\d=\underset{1\leq i\leq r}{\sum} d_i$.
Theorem\ref{th0-2}(1) is an immediate consequence of
Theorem\ref{th5-1} and Remark\ref{rem5-1} below.
\vskip 6pt

\begin{theorem}\label{th5-1}
\it Assume that $k$ is finitely generated over $\mathbb Q$ and
that $\d \geq n+1+\ccXZS$. Let $\kb$ be an algebraic closure of $k$
and put $\Sb=S\times_k\kb$.
Let $\etab$ be a geometric generic point of $\Sb$ and
let $U_{\etab}\subset X_{\etab}$ be the fibers of $\cU\subset\cX$ over $\etab$.
\begin{description}
\item[(1)]
Assuming $s\geq 1$,
$H^m_{et}(U_{\etab},\Ql(m))^{\pi_1(\Sb,\etab)}=
H^m_{et}(U_{\etab},\Ql(m))_{triv}.$
\item[(2)]
$H^m_{et}(X_{\etab},\Ql(m))_{prim}^{\pi_1(\Sb,\etab)}=0.$
\end{description}
\end{theorem}
\rm\medbreak

\begin{remark}\label{rem5-1}
\it Let $\wtd S \to S$ be a dominant morphism and
$\wtd{\cX}\supset \wtd{\cZ}$ be the base change. Then $\ccXZS=\ccXZtS$.
Hence Theorem\ref{th5-1} holds after replacing $\cX\supset\cZ$ by the base change.
In particular Theorem\ref{th5-1} holds if one replaces $\pi_1(\Sb,\etab)$ by any open
subgroup of finite index.
\end{remark}
\rm\medbreak

\begin{theorem}\label{th5-2}
\it Assume that $k=\mathbb C$ and $\d \geq n+1+\ccXZS$.
Let $U\subset X$ be the fibers of $\cU\subset\cX$ over a fixed base point
$0\in S(\mathbb C)$.
\begin{description}
\item[(1)]
Assuming $s\geq 1$,
$H^m(U,\mathbb Q(m))^{\pi_1(S,0)}=H^m(U,\mathbb Q(m))_{triv}.$
\item[(2)]
$H^m(X,\mathbb Q(m))_{prim}^{\pi_1(S,0)}=0.$
\end{description}
\end{theorem}
\rm\medbreak

First we deduce Theorem\ref{th5-1} from Theorem\ref{th5-2}.
By the Lefschetz principle we may assume that $k$ is a subfield of $\mathbb C$.
We fix an embedding $\kb(\etab)\hookrightarrow \mathbb C$ and let
$0\in \Sb(\mathbb C)$ be the corresponding
$\mathbb C$-valued point of $\Sb$. Write $\SC=S\otimes_k \mathbb C$ and put
$\UC=\cX\times_S \Spec(\mathbb C)$ via $\Spec(\mathbb C)\rmapo{0} \Sb\to S$.
We have the comparison isomorphisms ([SGA4, XVI Theorem4.1])
$$ H^m(\UC,\mathbb Q)\otimes \Ql \isom H^m_{et}(\UC,\Ql)\isom
 H^m_{et}(U_{\etab},\Ql)$$
that are equivariant with respect to the maps of topological and algebraic
fundamental groups:
$$ \pitop(\SC,0) \to \pialg(\SC,0) \to \pialg(\Sb,\etab).$$
The desired assertion follows at once from this.
\medbreak

In order to show Theorem\ref{th5-2} we need a preliminary. We assume $k=\mathbb C$. Let
$\San$ be the analytic site on $S(\mathbb C)$.
For a coherent sheaf $\cF$ on $\Szar$ let $\cF^{an}$ be the associated analytic
sheaf on $\San$. We introduce local systems on $\San$
$$\HQcU q p=R^q g_*\mathbb Q(p) \qaq \HCcU q =R^q g_*\mathbb C,$$
where $g:\cU\to S$ is the natural morphism.
Let $\HOcU q$ be the sheaf of holomorphic sections of $\HCcU q$
and let $F^p \HOcU q\subset \HOcU q$ be the holomorphic subbundle given
by the Hodge filtration on the cohomology of fibers of $\cU/S$.
We have the analytic Gauss-Manin connection
$$ \nabla\scs \HOcU q \to \WSan 1 \otimes \HOcU q \qwith
\Ker(\nabla)=\HCcU q$$
that satisfies
$\nabla(F^p\HOcU q) \subset \WSan 1 \otimes F^{p-1}\HOcU q.$
The induced map
$$ F^p\HOcU {p+q}/F^{p+1} \to \WSan 1 \otimes F^{p-1}\HOcU {p+q}/F^p$$
is identified with $(\nab p q)^{an}$ via the identification
$F^p \HOcU {p+q}/F^{p+1} = \HcUan p q$ where
$$  \nab p q \scs \HcU p q \to \WS {1}\otO \HcU {p-1}{q+1}$$
is defined in \S3. Therefore Theorem\ref{th3-1} implies the following.

\begin{theorem}\label{th5-3}
\it Let the assumption be as in Theorem\ref{th5-2}.
\begin{description}
\item[(1)]
If $s\geq 1$, $F^1\HOcU m \cap \HCcU m$
is generated over $\mathbb C$ by the image of
$$ H^m(U,\mathbb Q(m))_{triv}\hookrightarrow H^m(U,\mathbb Q(m))^{\pi_1(S,0)}
\simeq \Gamma(\San, \HQcU m m ) \hookrightarrow \HCcU m.$$
\item[(2)]
If $s=0$, $F^1\HOcX m \cap \HCcXp m=0$, where $\HCcXp q$ is the primitive part
of $\HCcX q$.
\end{description}
\end{theorem}
\rm\medbreak

Now we start the proof of Theorem\ref{th5-2}.
First we show Theorem\ref{th5-2}(2).
Write $H=H^{m}(X,\mathbb Q)_{prim}^{\pi_1(S,0)}$. We have the natural
isomorphism
$$ H\otimes\mathbb C \simeq \Gamma(\San,\HCcXp {m}).$$
By [D1] $H$ is a sub-Hodge structure of $H^{m}(X,\mathbb Q)$ so that
we have the Hodge decomposition
$$H\otimes\mathbb C=\underset{p+q=m}{\bigoplus} H^{p,q} \qwith
H^{p,q}\subset \Gamma(\San,F^p\HOcX {m} \cap \HCcXp {m}).$$
By Theorem\ref{th5-3} this implies $H^{p,q}=0$ for $p\geq 1$ which implies $H=0$
by the Hodge symmetry.

Next we show Theorem\ref{th5-2}(1). By Lemma\ref{lem5-1} below we have
$$ H^m(U,\mathbb Q(m))^{\pi_1(S,0)} \simeq \Gamma(\San,\HQcU m m)
\subset \Gamma(\San,\HCcU m \cap F^1\HOcU m).$$
Hence the assertion follows from Theorem\ref{th5-3}(1).
\medbreak

\begin{lemma}\label{lem5-1}
\it Let the assumption be as in Theorem\ref{th5-2}. Then
$H:=H^m(U,\mathbb Q)^{\pi_1(S,0)}$ is a submixed Hodge structure of
$H^m(U,\mathbb Q)$ and $F^1\HC=\HC$ where $\HC=H\otimes\mathbb C$.
\end{lemma}

\rm\begin{proof}
The fact that $H$ is a submixed Hodge structure follows from the theory
of mixed Hodge modules [SaM]. We show the second assertion.
We recall that the graded subquotients of the weight filtration of
$H^m(U,\mathbb Q)$ is given by
$$ Gr^W_{m+p}H^m(U,\mathbb Q)=
\underset{1\leq j_1<\cdots<j_p\leq s}{\bigoplus}
H^{m-p}(Z_{j_1}\cap\cdots\cap Z_{j_p},\mathbb Q(-p))_{prim}.
\quad (0\leq p\leq m)$$
Therefore it suffices to note that
$ H^{m}(X,\mathbb Q)_{prim}^{\pi_1(S,0)}=0$ by Theorem\ref{th5-2}(2).
\end{proof}
\medbreak

\def\tpi{\widetilde{\pi}}
\def\Skb{S_{\kb}}

Now we show Theorem\ref{th0-2}(2). It will follows from Theorem\ref{th5-4} below by using the
theory of Hilbert set (cf. [La]).

\begin{theorem}\label{th5-4}
\it Under the assumption of Theorem\ref{th0-2}(2), there exists an irreducible
variety $\tS$ over $k$ with a finite etale covering $\phi:\tS \to S$ for which
the following holds:
Let $\tpi=\pi\circ\phi$ and let $H\subset \mathbb P^N(k)$ be the subset of such
points that $\tpi^{-1}(y)$ is irreducible. For $\forall x\in S$ such that
$\pi(x)\in H$ and for any subgroup $G \subset Gal(k(\xb)/k(x))$
of finite index we have
$$H^m_{et}(\Uxb,\Ql(m))^G=H^m_{et}(\Uxb,\Ql(m))_{triv},$$
where $\xb$ is a geometric point over $x$ and $\Uxb$ is the fiber of $\cU/S$
over $\xb$.
\end{theorem}

\rm\begin{proof} (cf. [T] and [BE])
By choosing a $k$-rational point $0$ of $S$, we get the decomposition
$$\pi_1(S,\etab)=\pi_1(\Sb,\etab)\times \Gal(\kb/k),$$
where the notation is as in Theorem\ref{th5-1} and
$\pi_1(\Sb,\etab)$ is identified with the quotient of $\pi_1(S,\etab)$
which classifies the finite etale coverings of $S$
that completely decompose over $0$. Let
$$\Gamma=\Image\big(\pi_1(\Sb,\etab) \to
\GL_{\Ql}(H^m_{et}(U_{\etab},\Ql(m))\big).$$
From the fact that $\Gamma$ contains an $\ell$-adic Lie group as a subgroup
of finite index, we have the following fact (cf. [T]):
There exists a subgroup $\Gamma'\subset \Gamma$ of finite index such that
a continuous homomorphism of pro-finite group $P \to \Gamma$ is surjective
if and only if $P\to\Gamma\to \Gamma/\Gamma'$ is surjective as a map of sets.
Let $\phi:\tS\to S$ be a finite etale covering that corresponds to the inverse
image of $\Gamma'$ in $\pi_1(\Sb,\etab)$ and let $H$ be defined as in
Theorem\ref{th5-4}. Fix $x\in S$ with $\pi(x)\in H$ and let $\xb$ be a geometric
point of $x$. By choosing a ``path" $\xb\to \etab$, we get the map
$\Gal(k(\xb)/k(x))=\pi_1(x,\xb)\rmapo{\iota} \pi_1(S,\etab)$.
By the definition of $H$ the composite of $\iota$ with
$\pi_1(S,\etab)\to\Gamma\to\Gamma/\Gamma'$ is surjective so that
$\Gal(k(\xb)/k(x))$ surjects onto $\Gamma$ by the above fact.
It implies that for any subgroup $G \subset \Gal(k(\xb)/k(x))$ of finite index,
its image $\pi$ in $\Gamma$ is of finite index and hence that we have
$$ H^m_{et}(U_{\xb},\Ql(m))^G\isom
H^m_{et}(U_{\etab},\Ql(m))^\pi = H^m_{et}(U_{\etab},\Ql(m))_{triv}$$
where the last equality follows from Theorem\ref{th5-1} and Remark\ref{rem5-1}.
This completes the proof of Theorem\ref{th5-4}.
\end{proof}
\medbreak

Now Theorem\ref{th0-1}(2) is a consequence of Theorem\ref{th5-4} and the following result
(cf. [La]):

\begin{lemma}\label{lem5-2}
\it Let $k$ be a number field and let $V$ be an irreducible
variety over $k$. Let $\pi:V\to \mathbb P_k^N$ be an etale morphism
and let $H\subset \mathbb P_k^N(k)$ be the subset of such points $x$
that $\pi^{-1}(x)$ is irreducible.
Let $\Sigma$ be any finite set of primes of $k$ and let $k_v$ be the
completion of $k$ at $v$. Then the image of $H$ in
$\prod_{v\in \Sigma} \mathbb P_k^N(k_v)$ is dense.
\end{lemma}
\rm

\vskip 20pt

\setcounter{section}{6}\setcounter{subsection}{0}
\begin{center}{\bf\arabic{section}.
Implication on injectivity of regulator maps for $K_1$ of
surfaces}\end{center}
\medbreak

Let $X$ be a projective smooth surface over a field $k$.
Let $U\subset X$ be the complement of a simple normal crossing
divisor $Z\subset X$. In this section we discuss an implication of
the surjectivity of $\regU 2$ and $\regetU 2$ on the regulator maps
for $CH^2(X,1)$. Recall that $CH^2(X,1)$ is by definition the cohomology of
the following complex
$$ K_2(k(X)) \rmapo{\dt} \bigoplus_{C\subset X} k(C)^*
\rmapo{\dd} \bigoplus_{x\in X} \mathbb Z,$$
where the sum on the middle term ranges over all irreducible curves on $X$
and that on the right hand side over all closed points of $X$.
The map $\dt$ is the so-called tame symbol and $\dd$ is the sum of divisors
of rational functions on curves. Thus an element of $\Ker(\dd)$ is given by a
finite sum
$\sum_{i} (C_i,f_i)$, where $f_i$ is a non-zero rational function on
an irreducible curve $C_i\subset X$ such that $\sum_{i} div(f_i)=0$ on $X$.
We recall that $K_2(k(X))$ is generated as an abelian groups by symbols
$\{f,g\}$ for non-zero rational functions $f,g$ on $X$ and that
$$ \dt(\{f,g\})=((f)_0,g)+((f)_{\infty},1/g)+((g)_0,1/f)+((g)_{\infty},f),$$
where $(f)_0$ (resp. $(f)_{\infty}$) is the zero (resp. pole) divisor of $f$.
An important tool to study $CH^2(X,1)$ is the regulator map:
In case $k=\mathbb C$ it is given by
$$ \regD X \scs CH^2(X,1) \to H^3_D(X,\mathbb Z(2)),$$
where the group on the right hand side is the Deligne cohomology group
(cf. [EV] and [J2]).
Assuming the first Betti number $b_1(X)=0$, we have the following explicit
desciption of $\regD X$. Take $\alpha=\sum_{i} (C_i,f_i)\in \Ker(\dd)$.
Under the isomorphism (cf. [EV, 2.10])
$$ H^3_D(X,\mathbb Z(2))\simeq \frac{H^2(X,\mathbb C)}{H^2(X,\mathbb Z(2))+ F^2H^2(X,\mathbb C)}
\simeq \frac{F^1H^2(X,\mathbb C)^*}{H_2(X,\mathbb Z)},$$
$\regD X(\alpha)$ is identified with a linear function on complex valued
$C^{\infty}$-forms $\omega$ and we have
$$ \regD X(\alpha)(\omega)=
\frac{1}{2\pi\sqrt{-1}} \sum_i \int_{C_i-\gamma_i} \log(f_i)\omega
+ \int_{\Gamma} \omega,$$
where $\gamma_i:=f_i^{-1}(\gamma_0)$ with $\gamma_0$, a path on
$\mathbb P^1_{\mathbb C}$ connecting $0$ with $\infty$ and $\Gamma$ is a real
piecewise smooth $2$-chain on $X$ such that
$\partial\Gamma =\cup_i \gamma_i$ which exists due to the assumption
$\alpha\in \Ker(\dd)$ and $b_1(X)=0$.
\vskip 5pt
In case $k$ is finitely generated over $\mathbb Q$ we have the regulator map
$$\regcont X \scs CH^2(X,1)\otimes\Zl\to
H_{cont}^3(X,\Zl(2)),$$
where $H^i_{cont}$ denotes the continuous etale cohomology of $X$ (cf. [J3]).
\medbreak

Now let $Z\subset X\supset U$ be as in the begining of this section and write
$Z=\cup_{1\leq i\leq r} Z_i$ with $Z_i$, smooth irreducible curves intersecting
transversally with each other.
We consider the following subgroup of $\Ker(\dd)$
$$ CH^1(Z,1)=\Ker(\bigoplus_{1\leq i\leq r} k(Z_i)^* \rmapo{\dd}
\bigoplus_{x\in Z} \mathbb Z).$$
By the localization theory for higher Chow group we have the exact sequence
$$ CH^2(U,2) \rmapo{\dz} CH^1(Z,1) \rmapo{\iota} CH^2(X,1)$$
where the first map coincides up to sign with the composite of
the natural map $CH^2(U,2) \to K_2(k(X))$ and $\dt$.

\begin{theorem}\label{th6-1} \it
\begin{description}
\item[(1)]
Assume $k=\mathbb C$ and that there exists a subspace
$\Delta\subset CH^2(U,2)\otimes \mathbb Q$ such that the restriction of
$\regU 2$ on $\Delta$ is surjective. Let $\alpha\in CH^1(Z,1)$ and assume
$\regD X(\iota(\alpha))=0$. Then $\alpha\in \dz(\Delta)$ in
$CH^1(Z,1)\otimes\mathbb Q$. In particular
$\iota(\alpha)=0\in CH^2(X,1)\otimes\mathbb Q$.
\item[(2)]
Assume that $k$ is finitely generated over $\mathbb Q$.
The analogous fact holds for $\regetU 2$ and $\regcont X$.
\end{description}
\end{theorem}
\rm\medbreak

The following result is a direct consequence of Theorem\ref{th6-1} and
Definition\ref{def0-1}.
Let $\cZ\subset \cX$ be as (0-4) in the introduction and let $\Zt\subset \Xt$
be the fibers over $t\in S$.

\begin{corollary}\label{cor6-2}
\it Let the assumption be as in Theorem\ref{th0-1}.
For any $t\in S\setminus \SNL$, $\regD {\Xt}$ restricted on the image of
$CH^1(\Zt,1)$ is injective modulo torsion.
\end{corollary}

\begin{remark}\label{rem6-1}
In the forthcoming paper [AS2] we study the Noether-Lefschetz
locus in case (0-4) is the universal family of
smooth surfaces of degree $d$ and three hyperplanes in $\mathbb P^3$
intersecting transversally. In this case it is shown that there exist
$t\in S\setminus \SNL$ such that the image of $CH^1(\Zt,1)$ in
$CH^2(\Xt,1)$ is non-torsion if $d$ is large enough. Thus Corollary\ref{cor6-2} has indeed non-trivial implication on the injectivity of the regulator map.
\end{remark}
\medbreak

Now we start the proof of Theorem\ref{th6-1}.
The idea of the following proof is taken from [J1, 9.8].
We only treat the first assertion. The second is proven by the same way.
We have the commutative diagram (cf. [Bl] and [J2, 3.3 and 1.15]).
$$\begin{matrix}
CH^2(X,2) & \to & CH^2(U,2) &\to& CH^1(Z,1) &\to& CH^2(X,1) \\
\downarrow\rlap{$\chD 22X$} && \downarrow\rlap{$\chD 22U$}&&
\downarrow\rlap{$\chD 11Z$} && \downarrow\rlap{$\chD 21X$}\\
H^2_D(X,\mathbb Z(2)) &\rmapo{\iota}& H^2_D(U,\mathbb Z(2)) &\to&
H^3_{D,Z}(X,\mathbb Z(2)) &\to& H^3_D(X,\mathbb Z(2)) \\
\end{matrix}$$
Here the horizontal sequences are the localization sequences
for higher Chow groups and Deligne cohomology groups and
they are exact. The vertical maps are regulator maps.
By a simple diagram chasing it suffices to show that $\chD 11Z$ is
injective and
$\chD 22U(\Delta)+\Image(\iota)$ spans $H^2_D(U,\mathbb Q(2))$.
In order to show the first assertion we note the commutative diagram
$$\begin{matrix}
0\to&\underset{1\leq i\leq r}{\bigoplus} CH^1(Z_i,1) &\to& CH^1(Z,1) &\to&
\underset{x\in W}{\bigoplus} \mathbb Z\\
&\downarrow\rlap{$\chD 11{Z_i}$} && \downarrow\rlap{$\chD 11Z$}
&&  \llap{$\simeq$}\downarrow\rlap{$\chD 00x$}\\
0\to&\underset{1\leq i\leq r}{\bigoplus} H^1_{D}(Z_i,\mathbb Z(1)) &\to&
H^3_{D,Z}(X,\mathbb Z(2)) &\to&
\underset{x\in W}{\bigoplus} H^0_{D}(x,\mathbb Z(0))\\
&\downarrow\rlap{$\simeq$}&&\downarrow\rlap{$=$}&&\downarrow\rlap{$\simeq$}\\
&\underset{1\leq i\leq r}{\bigoplus} H^3_{D,Z_i}(X,\mathbb Z(2)) &\to&
 H^3_{D,Z}(X,\mathbb Z(2)) &\to&
\underset{x\in W}{\bigoplus}H^4_{D,x}(X,\mathbb Z(2))\\
\end{matrix}$$
where $W=\coprod_{1\leq i\not=j\leq r} Z_i\cap Z_j$.
The horizontal sequences come from the Mayer-Vietoris spectral sequence
and they are exact. Thus the desired assertion follows from the fact that
$\chD 11{Z_i}$ is an isomorphism (cf. [J2, 3.2]).
To show the second assertion we recall the exact sequence
(cf. [EV, 2.10])
$$ 0 \to H^1(U,\mathbb C)/H^1(U,\mathbb Z(2)) \rmapo{\gamma}
H^2_D(U,\mathbb Z(2))
\rmapo{\beta} H^2(U,\mathbb Z(2))\cap F^2H^2(U,\mathbb C) \to 0$$
and the same sequence with $U$ replaced by $X$. We have
$\chB 22 U=\beta\cdot \chD 22U$ and the commutative diagram
$$\begin{matrix}
\mathbb C^*\otimes CH^1(U,1) &\to& CH^1(U,1)\otimes CH^1(U,1)\\
\downarrow\rlap{$\log\otimes \chB 11 U$}&&\downarrow\rlap{\text{product}}\\
\mathbb C/\mathbb Z(1)\otimes H^1(U,\mathbb Z(1)) && CH^2(U,2)\\
\downarrow\rlap{$\simeq$}&&\downarrow\rlap{$\chD 22U$}\\
H^1(U,\mathbb C)/H^1(U,\mathbb Z(2)) &\rmapo{\gamma}& H^2_D(U,\mathbb Z(2))\\
\end{matrix}$$
Therefore the desired assertion follows from the surjectivity of
$$CH^1(U,1) \rmapo{\chB 11 U} H^1(U,\mathbb Z(1))/H^1(X,\mathbb Z(1))$$
which is easily seen.

\end{document}